\documentclass[11pt, a4paper]{amsart}
\usepackage{amssymb}
\usepackage{amsmath}
\usepackage{amsthm}
\usepackage[fleqn]{mathtools}
\usepackage{amsfonts}
\usepackage{stmaryrd}
\usepackage[english]{babel}
\usepackage[utf8]{inputenc}
\usepackage{enumerate}
\usepackage{graphicx,pinlabel,cases}
\usepackage{float}

\newtheorem{theorem*}{Theorem}
\newtheorem{corollary*}{Corollary}
 \theoremstyle{definition}
 \theoremstyle{remark}

 \numberwithin{equation}{section}
 
\begin{document}

\title{Estimating the Cheeger constant using machine learning}

\author[A. Jain]{Ambar Jain}
\address{Department of Physics\\
Indian Institute of Science Education and Research Bhopal\\
Bhopal Bypass Road, Bhauri \\
Bhopal 462 066, Madhya Pradesh\\
India}
\email{ambarj@iiserb.ac.in}
\urladdr{http://home.iiserb.ac.in/$_{\widetilde{\phantom{n}}}$ambarj/}

\author[S. Pal]{Shivam Pal}
\address{Department of Physics\\
Indian Institute of Science Education and Research Bhopal\\
Bhopal Bypass Road, Bhauri \\
Bhopal 462 066, Madhya Pradesh\\
India}
\email{cvampal@iiserb.ac.in}

\author[K. Rajeevsarathy]{Kashyap Rajeevsarathy}
\address{Department of Mathematics\\
Indian Institute of Science Education and Research Bhopal\\
Bhopal Bypass Road, Bhauri \\
Bhopal 462 066, Madhya Pradesh\\
India}
\email{kashyap@iiserb.ac.in}
\urladdr{https://home.iiserb.ac.in/$_{\widetilde{\phantom{n}}}$kashyap/}



\begin{abstract}
In this paper, we use machine learning to show that the Cheeger constant of a connected regular graph has a predominant linear dependence on the largest two eigenvalues of the graph spectrum. We also show that a trained deep neural network on graphs of smaller sizes can be used as an effective estimator in estimating the Cheeger constant of larger graphs.
\end{abstract}

\maketitle

\section{Introduction}
 Let $G = (V,E)$ be a finite, simple, connected and undirected $k$-regular graph with $|G| = n$. It is a well known fact from basic algebraic graph theory~\cite{B2,G1} that the eigenvalues $\lambda_i(G)$, $0 \leq i \leq n-1$, of the adjacency matrix $A(G)$ of $G$ are real and can be ordered as:
$$k= \lambda_0(G) \geq \lambda_1(G) \geq \dotsc \geq \lambda_{n-1}(G) \geq -k.$$
For each $F \subset V$, let $\partial F := \{\{ u, v \} \in E(X) : u \in F, v \in V \setminus F\}$. Then
the number 
\begin{equation}
\label{eq:bound1}
h(G) = \min_{F \subset V, \,|F| \leq \frac{|V|}{2}} \frac{|\partial F|}{|F|},
\end{equation} 
is called the \textit{Cheeger constant} (or the \textit{isoperimetric constant} or the \textit{edge expansion constant}) of the graph $G$. The Cheeger constant is a measure of the connectivity of the graph $G$. Families of regular graphs with Cheeger constants bounded below by a positive constant also known as \textit{expander families} have been widely studied (see~\cite{HLW,LUB,AL,RM1} and the references therein) due to their applications to communication networks. 

The computation of $h(G)$ of an arbitrary finite graph is a well known~\cite{GJS,LR,BM} NP hard problem. However, for a $k$-regular graph $G$ of size $n$, we use machine learning to answer the following natural questions. 

\begin{enumerate}[(a)]
\item Does the dependence of $h(G)$ on $\lambda_0(G) = k$ and $\lambda_1(G)$ stronger than what the known bounds indicate?
\item Is this dependence predominantly linear or non-linear? 
\item Is there a strong dependence of $h(G)$ on $\lambda_i(G)$, for $2 \leq i \leq n-1$?
\item Can these dependencies be used to estimate $h(G)$ for large $n$ with greater efficiency? 
\end{enumerate}

We begin by providing data which shows that in general these known bounds for $h(G)$ deviate significantly from its actual value. By considering random regular graphs of sizes 12 through 30, we apply machine learning via deep neural networks and linear regression to make the following statistical observations: 
\begin{enumerate}[(i)]
\item $h(G)$ has predominant linear dependence on $\lambda_0(G)$ and $\lambda_1(G)$. Moreover, as $|G|$ increases, this dependence appears to approach the linear function $\frac{1}{2}\lambda_0(G) - \frac{1}{3}\lambda_1(G)$. This linearity is more pronounced when the spectral gap is large.
\item Its dependence on $\lambda_i(G)$, for $2 \leq i \leq n-1$ is insignificant.
\item We demonstrate that a trained deep neural network on graphs of smaller sizes can be used as an effective estimator for Cheeger constants of larger graphs where computation times using classical algorithms are large.
\end{enumerate}

The paper is organized as follows. In Section~\ref{sec:bound_analysis}, we analyze whether some well known bounds can be used as effective estimators for $h(G)$. In Section~\ref{sec:linearpredict}, we determine whether the dependence of $h(G)$ on $\lambda_0(G)$ and $\lambda_1(G)$  is predominantly linear. In Section~\ref{sec:non_linear_dep}, we use machine learning to examine whether $h(G)$ has a nonlinear dependence on $\lambda_0(G$)  and $\lambda_1(G)$, and also study its relation to $\lambda_i(G)$, for $2 \leq i \leq n-1$. Finally, in Section~\ref{sec:pred_Cheeger_constant}, we explore whether deep neural networks trained on graphs of smaller sizes can be used as viable estimators for Cheeger constants of larger graphs.

\section{Numerical analysis of known bounds}
\label{sec:bound_analysis}

We consider a dataset of random regular graphs of sizes 12 through 30 for our analysis. This dataset was generated by using a Python package that implements the algorithm described in~\cite{SW}. The number of graphs considered for $n=12$ was limited by the total number of available graphs, while for $n>20$, the limitation came from long computation time for $h(G)$. In all other cases, we have considered at least $20,000$ random graphs of varying regularity. The number of graphs considered for analysis for each $n$ is shown in the second column of Table \ref{tab:bound-analysis}.

The Cheeger is related to the \textit{spectral gap} $k - \lambda_1(G)$ of a $k$-regular graph $G$ by the following inequality (see~\cite[Proposition 1.84]{KS}): 
\begin{equation}
\label{eq:bound2}
\frac{k - \lambda_1(G)}{2} \leq h(G) \leq \sqrt{2k(k-\lambda_1(G))}.
\end{equation}  Mohar~\cite{BM} showed that 
\begin{equation}
\label{eq:bound3}
h(G) \leq \begin{cases}
					\displaystyle 
 					\frac{k}{2} \left[ \frac{n}{n-1} \right], & \text{if } n \text{ is even, and} \\ \\
 					\displaystyle 
 					\frac{k}{2}\left[\frac{n+1}{n-1}\right] & \text{if } n \text{ is odd.}
				    \end{cases},
 \end{equation}
and when $G \neq K_1,\, K_2, \text{ or } K_3$ (where $K_i$ denotes the complete graph on $i$ vertices), he showed that
\begin{equation}
\label{eq:bound4}
h(G) \leq \sqrt{k^2 - \lambda_1(G)^2}.
\end{equation}

For each graph $G$, we compute the lower bound on $h(G)$ as given by Eqn.~(\ref{eq:bound2}), and an upper bound, which is the lowest of the upper bounds appearing in (\ref{eq:bound2})-(\ref{eq:bound4}).  For each of these estimators, we calculate its deviation $\Delta h$ from the true value of $h(G)$ as given in the equation below:
\begin{equation}
\Delta h_{\rm est.} = \left \vert \frac{h_{\rm est.} - h(G)}{h(G)} \right \vert \, ,
\end{equation}
where, $h_{\rm est.}$ refers to the estimator of $h(G)$. For the analysis in this section, $h_{\rm est.}$ corresponds to either the upper bound or the lower bound. The mean values of $\Delta h_{\rm est.}$ (which we denote by $\left \langle \Delta h_{\rm lower} \right \rangle$ and $\left \langle\Delta h_{\rm upper}\right \rangle$ respectively)
 for each $n$ is shown in Table \ref{tab:bound-analysis} below. 

\begin{table}[h]
\centering
\begin{tabular}{|c|c|c|c|}
\hline
 $n$ &  \# of Graphs & $\left \langle \Delta h_{\rm lower} \right \rangle$ &  $\left \langle\Delta h_{\rm upper}\right \rangle$\\
\hline
 12 & 15176 & 0.18 & 0.61 \\
 13 & 55128 & 0.23 & 0.61 \\
 14 & 115663 & 0.18 & 0.60 \\
 15 & 118702 & 0.22 & 0.63 \\
 16 & 22635 & 0.18 & 0.65 \\
 17 & 20024 & 0.21 & 0.61 \\
 18 & 35774 & 0.18 & 0.66 \\
 19 & 20436 & 0.20 & 0.59 \\
 20 & 56016 & 0.18 & 0.64 \\
 21 & 1606 & 0.19 & 0.56 \\
 22 & 1626 & 0.17 & 0.61 \\
 23 & 1636 & 0.19 & 0.55 \\
 24 & 1825 & 0.16 & 0.59 \\
 25 & 1385 & 0.18 & 0.53 \\
 26 & 1829 & 0.16 & 0.57 \\
 27 & 1722 & 0.17 & 0.52 \\
 28 & 1097 & 0.16 & 0.57 \\
 29 & 958 & 0.20 & 0.62 \\
 30 & 872 & 0.16 & 0.59 \\
\hline
\end{tabular}
\caption{Graph data considered in the analysis of this paper and the average deviation in bounds: The second column shows the number of graphs considered in the analysis in this paper for each $n$. For $n \le 20$, at least 20,000 graphs were considered for each $n$ with exception to $n=12$, where the total number of available graphs is less than 20,000. For $n>21$, we tried to accumulate at least about 1000 graphs with the exceptions of $n=29$ and $n=30$.}
\label{tab:bound-analysis}
\end{table}
\noindent We note that, on an average, the lower bound deviates from the true value of $h(G)$ by about $20\%$, while the upper bound deviates at about $60\%$. This deviation marginally reduces for large values of $n$. The table indicates that the bounds considered are not efficient estimators for $h(G)$. In the following section, we consider linear regression to construct a better estimator for $h(G)$.

\section{Linear regression analysis and prediction}
\label{sec:linearpredict}
\begin{figure}[htbp]
\centering
\includegraphics[width=65ex]{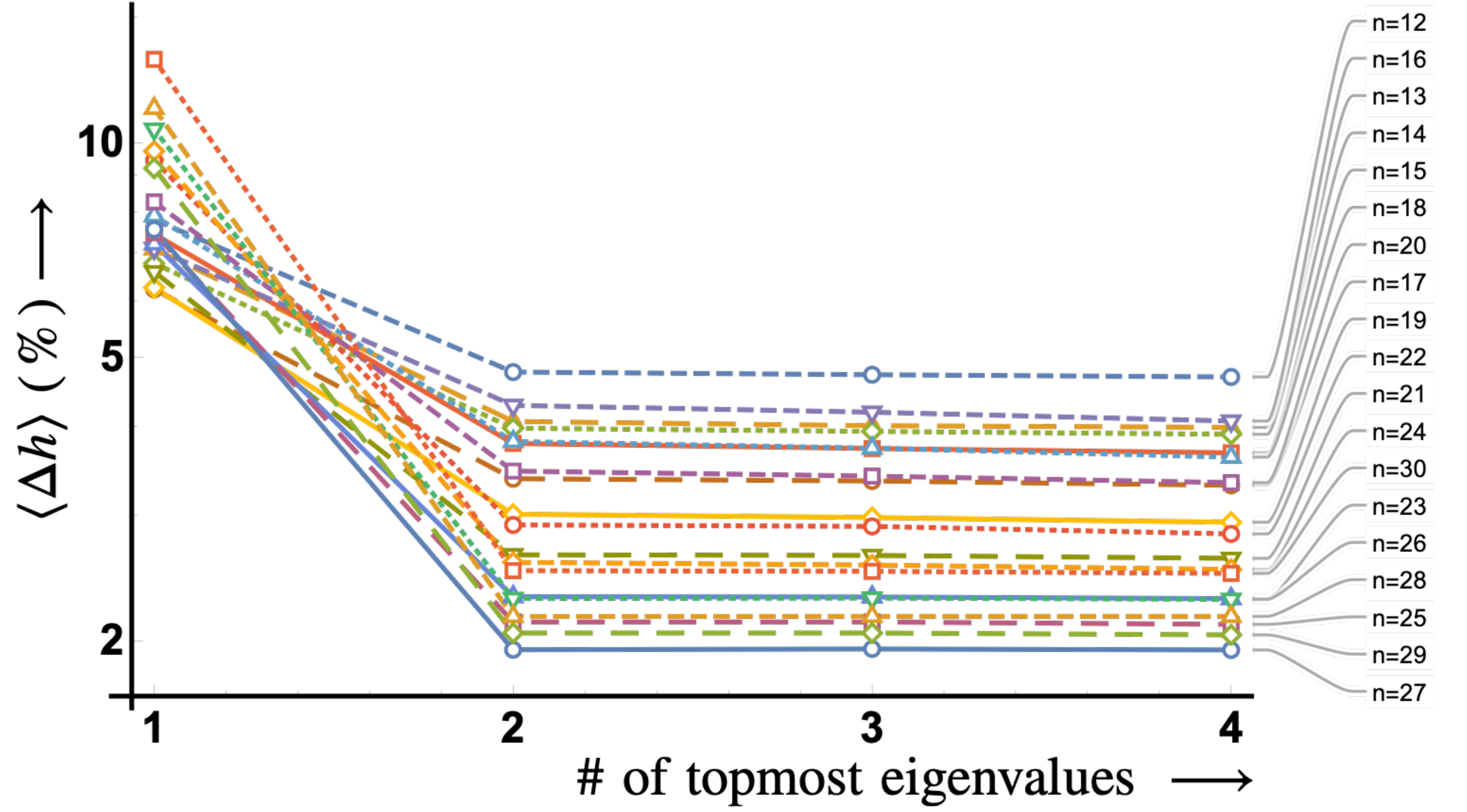}
\caption{Linear regression on Cheeger constant: Graph shows average deviation of $h(G)$ from the estimate obtained through a linear fit for topmost, top two, top three and top four eigenvalues. Points are joined by lines to guide the eye. Log scale is used on the $y$-axis to stretch the scale. There is no considerable improvement in the Cheeger estimate from linear regression beyond $\lambda_0(G)$ and $\lambda_1(G)$.}
\label{fig:linear-regression}
\end{figure}

In this section, we want to determine whether the relationship between $h(G)$ and $\lambda_0(G)$ and $\lambda_1(G)$ is predominantly linear. To begin with, we analyze whether $h(G)$ can be estimated reasonably well by a linear function of the largest $m$ eigenvalues, for $1 \leq m \leq 4$. For each $m$, we calculate the mean deviation $\langle \Delta h \rangle$, where we use the fitted linear regression function as the estimator. The results for this analysis are presented in Fig.~\ref{fig:linear-regression} below for various values of $n$.

 It is evident from the graph that adding the third and fourth eigenvalue to the analysis does not significantly reduce $\langle \Delta h\rangle$. This shows that a linear function of just the two largest eigenvalues estimates $h(G)$ fairly accurately Interestingly, the average deviation $\langle \Delta h\rangle$ reduces gradually with increase in $n$ coming down to about $2\%$ for $n \approx 30$. This observation confirms that the relationship between the two largest eigenvalues and $h(G)$ is mostly linear.

\begin{figure}[H]
\centering
\includegraphics[width=55ex]{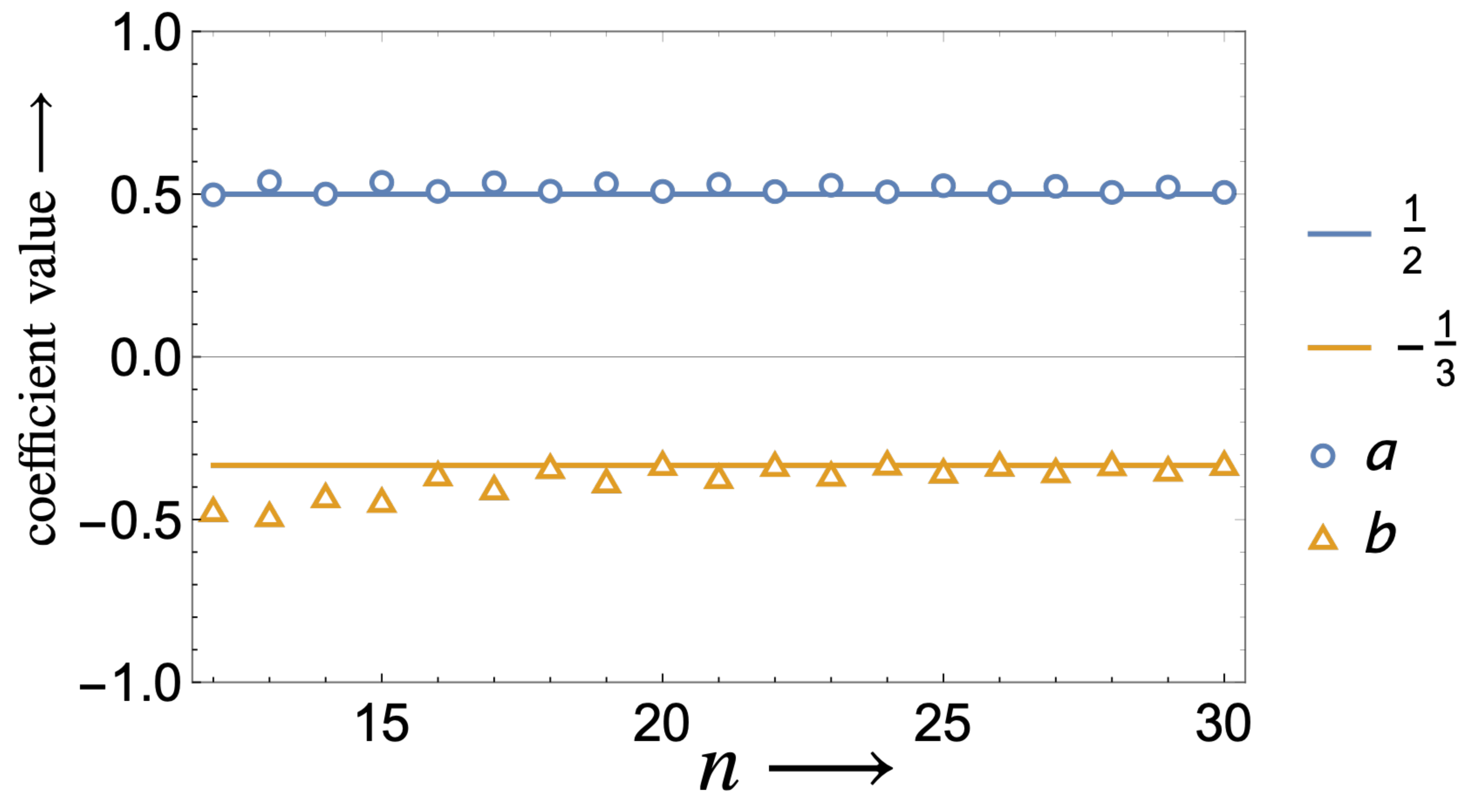}
\caption{Coefficients of linear regression on $h(G)$. }
\label{fig:linearcoeff}
\end{figure}

The regression coefficients of $\lambda_0(G)$ appears to converge to $\frac{1}{2}$ as $n$ increases, while the coefficient of $\lambda_1(G)$ appears to converge to $-\frac{1}{3}$. The coefficients $a$ and $b$ of the model $a \lambda_0(G) + b \lambda_1(G) + c$ are plotted in Fig.~\ref{fig:linearcoeff} below for each $n$ along with lines corresponding to $1/2$ and $-1/3$ for reference.

\noindent  This suggests a universality in the linear relationship, which is almost independent of $n$. This observation motivates us to test the linear model on $\lambda_0(G)$ and $\lambda_1(G)$ for the prediction of $h(G)$ for larger $n$, where its computation is challenging. We train the linear regression model on the available data for $ n=12,13,16,17$ and then use it to predict $h(G)$ for other $n$. Using the trained linear model as the estimator, we show the mean deviation  $\langle \Delta h \rangle$ in Fig.\ref{fig:linear-prediction} below. 
\begin{figure}[H]
\centering
\includegraphics[width=70ex]{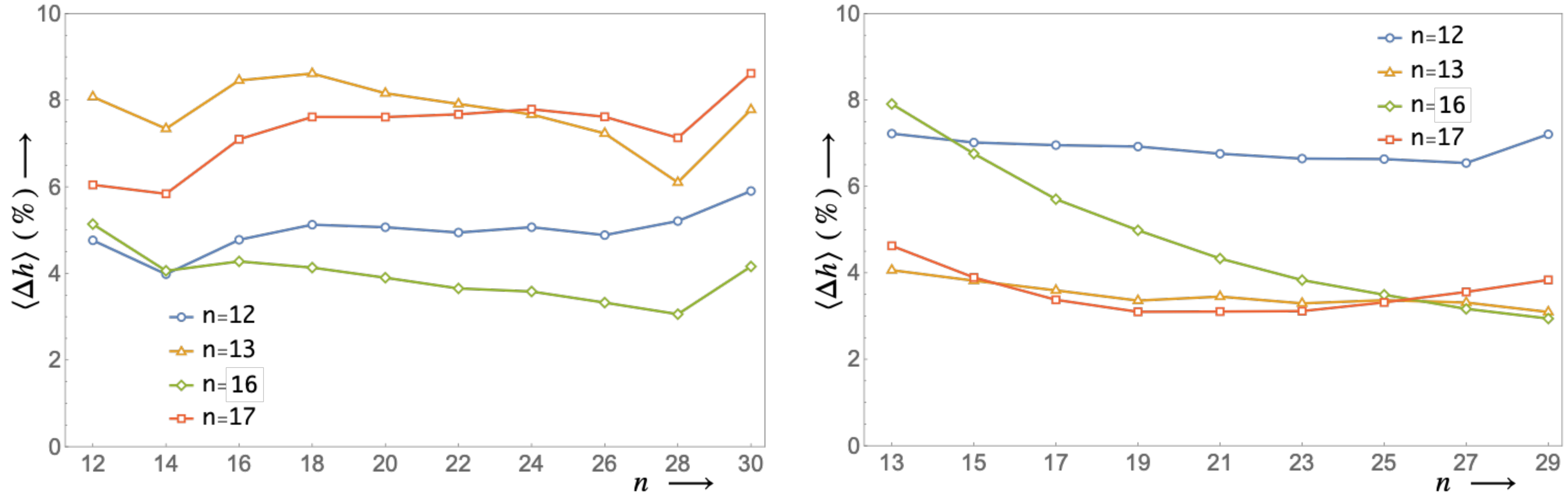}
\caption{Predicting with linear regression: Linear models trained on $h(G)$ data for $n=12,\, 13,\, 16,\, 17$ are used to predict Cheeger constants for graphs of other $n$. Average fractional deviation of the model from true value of $h(G)$ is shown for each $n$. Linear models trained on even (resp. odd) values of $n$ work better for the prediction of $h(G)$ for even (resp. odd) $n$.}
\label{fig:linear-prediction}
\end{figure}
\noindent The left panel shows prediction for even $n$, while the right panel shows prediction for odd $n$. We make the following observations
\begin{enumerate}
\item In general, for large $n$, linear regression with $\lambda_0(G)$ and $\lambda_1(G)$ appears to be a reasonable estimator for $h(G)$. 
\item The prediction is slightly more accurate when regression on odd $n$ (resp. even $n$) is used to predict the $h(G)$ for larger values of odd $n$ (resp. even $n$.)
\item Average deviation $\langle \Delta h \rangle$ is typically 4-5\% for odd-odd and even-even predictions for the entire range of $n$ considered.
\item{It also appears that $n=16$ and $n=17$ linear models are slightly better over $n=12$ and $n=13$ models respectively for even-sized and odd-sized graphs respectively.} This indicates that, for training a predictive model, we should opt for largest possible even and odd $n$ for which the Cheeger constant data is available.
\end{enumerate}

\section{Estimation of Cheeger constant using machine learning}
\label{sec:non_linear_dep}
\begin{figure}[H]
\centering
\includegraphics[width=65ex]{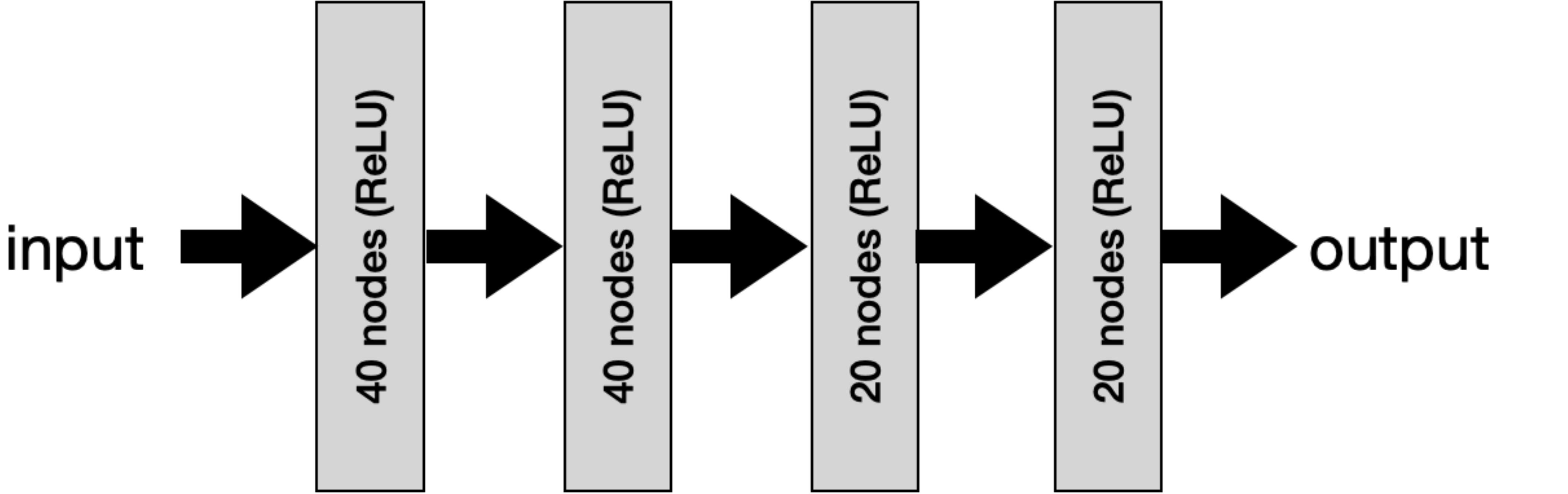}
\caption{Neural Network architecture used in this paper.}
\label{fig:neuralNet}
\end{figure}

In this section, we study the data on $h(G)$ using machine learning methods with deep neural networks, mainly to answer following two questions.
\begin{enumerate}
\item{Does $h(G)$ have a non-linear dependence on $\lambda_0(G)$ and $\lambda_1(G)$?}
\item{Does $h(G)$ has any significant dependence on other eigenvalues?}
\end{enumerate}
We expect that machine learning techniques will be able to identify non-linear dependencies that were not visible through linear regression. We randomly take 40\% of our dataset for $12 \leq n \leq 30$ and train a deep neural network shown in Fig.~\ref{fig:neuralNet} below \footnote{We have observed that other similar choices of neural net produce similar results presented in this section, as is the case with any machine learning problem. Several results in this paper can also be produced using a less deeper network. Our choice of neural network here works for all the results presented here.} using ADAM optimizer. 

\noindent The remaining $60\%$ of the dataset is used for validation. The trained neural net essentially provides an approximate non-linear map between the input eigenvalues and the expected Cheeger constant.
The validation ensures that there is no memorization done by the neural net and that it is truly capturing features of the data. 
Fig.~\ref{fig:DNN12} below shows training and validation histograms of $\Delta h$ for $n=12$, for both the cases of trainings done with the largest two and the largest four eigenvalues.
\begin{figure}[htbp]
\centering
\includegraphics[width=55ex]{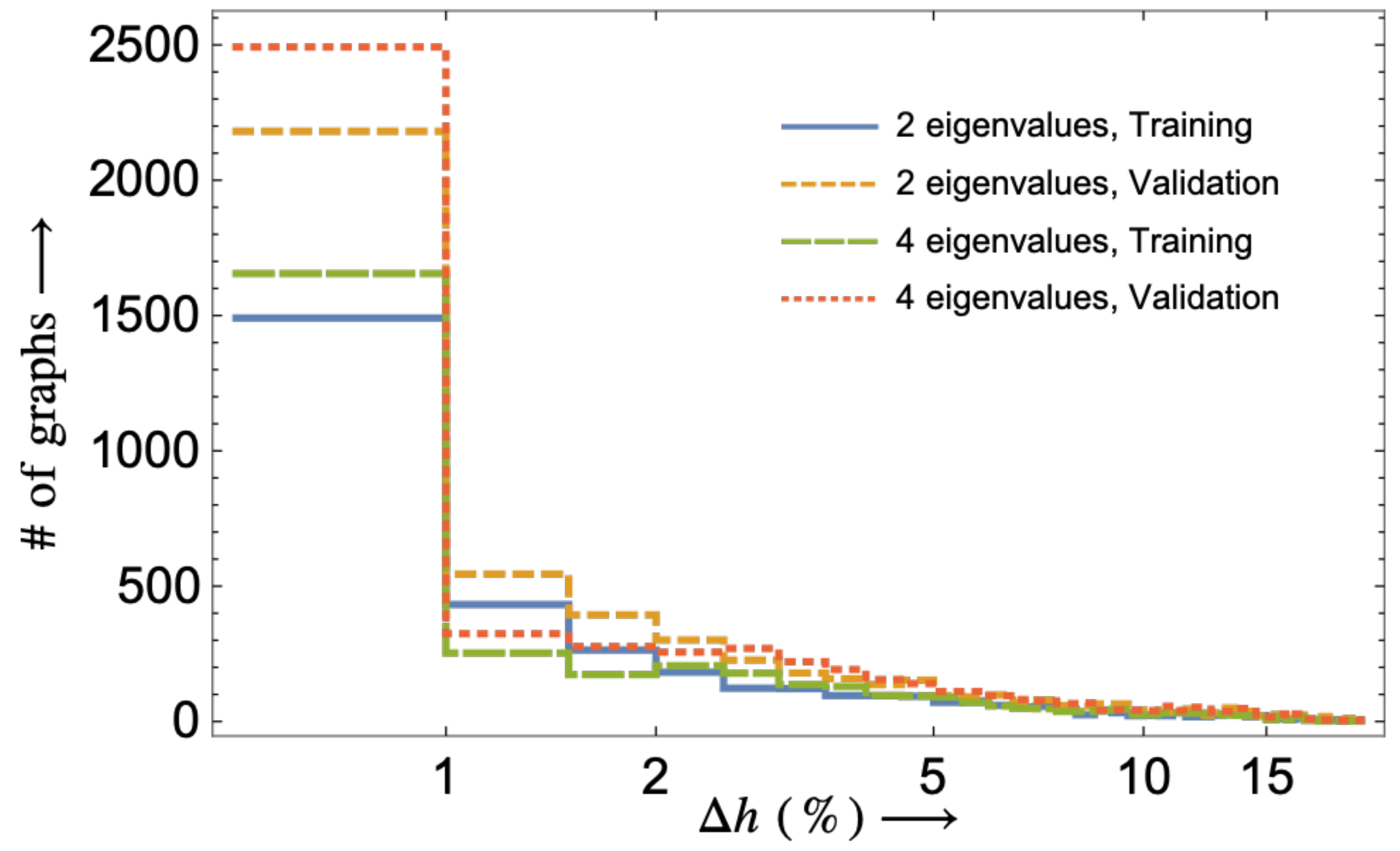}
\caption{Training and Validation Histograms for training deep neural network on $n=12$ graph data. Each bin size corresponds to 0.5\% of $\Delta h$. Mean deviation for both training and validation for both cases of $2$ and $4$ eigenvalues is about $2.5\%$. }
\label{fig:DNN12}
\end{figure}

We make the following observations:
\begin{enumerate}
\item $\lambda_0(G)$ and $\lambda_1(G)$ have a very strong correlation with $h(G)$. Furthermore, there appears to be a small non-linear dependence on $\lambda_0(G)$ and $\lambda_1(G)$ which accounts for about 2.5\% improvement over the linear regression. The average deviation $\langle \Delta h \rangle$ is about 2.5\% in both the training and validation data sets for deep neural net (DNN) model while it was about 5\% for the linear model.
\item We do not observe any significant improvement for the estimation of $h(G)$ while considering largest four eigenvalues over $\lambda_0(G)$ and $\lambda_1(G)$. In both cases $\langle \Delta h \rangle \approx 2.5\%$ with small fluctuations in each attempt of training. Using other subsets of the spectrum, including the full spectrum, does not seem to improve the training and validation errors beyond what are observed by considering just $\lambda_0(G)$ and $\lambda_1(G)$.
\item Similar exercise done for other graph sizes between $n=13$ to $n=20$ has similar results. Mean and standard deviation for $\Delta h$ for these cases is plotted in Fig.~\ref{fig:DNNfor13to20} below for both training and validation, reaffirming the observations made above.

\begin{figure}[H]
\includegraphics[width=55ex]{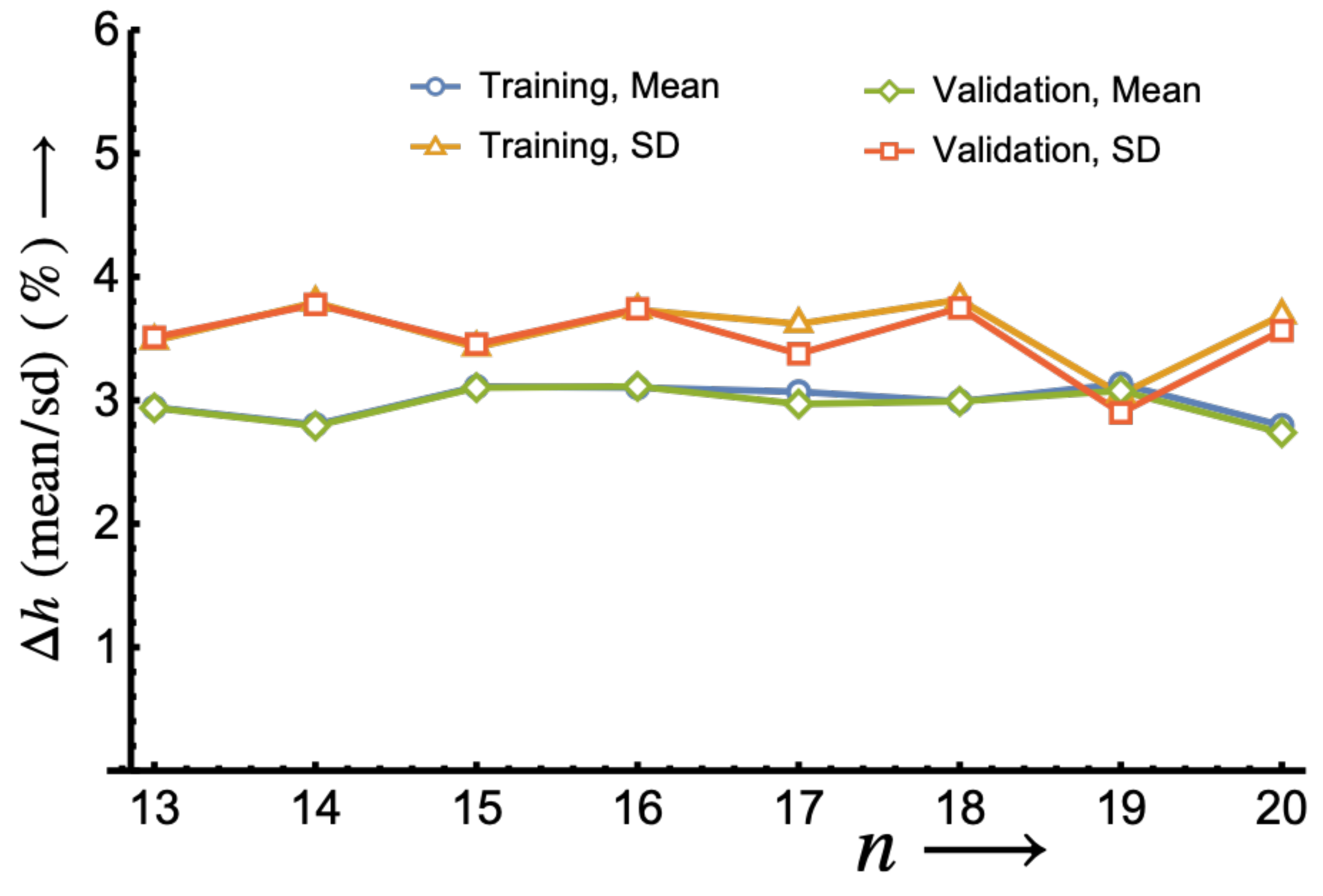}
\caption{Training and validation mean and their standard deviations for neural network model trained with $\lambda_0(G)$ and $\lambda_1(G)$ for $n=13$ to $20$. }
\label{fig:DNNfor13to20}
\end{figure}

\item Studying the trained deep neural network reveals that $h(G)$ has largely linear dependence on $\lambda_0$ and $\lambda_1$ when the spectral gap is large, while it exhibits non-linear dependence when the spectral gap is small. 
\end{enumerate}

\noindent We conclude that $\lambda_0(G)$ and $\lambda_1(G)$ suffice to estimate $hG)$ reliably. 

\section{Predicting $h(G)$ using Machine Learning}
\label{sec:pred_Cheeger_constant}
The most interesting application of this work is to predict Cheeger constant for large regular graphs, where it is computationally inefficient to calculate Cheeger constant but computationally efficient to calculate the spectrum. To achieve this, we train a neural net for small graphs where it is possible to calculate Cheeger constant in reasonable computation time. We then use this trained net to predict Cheeger constant for the large graph. We moderately train the deep neural network shown in the previous section for 50 epochs\footnote{The training was stopped after 50 epochs as compared to about 500 epochs (optimization stopping automatically when loss stops improving) done in the previous section. This ensures that the network learns the significance of top two eigenvalues and not the information about the $n$. Maximal training to about 500 epochs optimizes the network to estimate Cheeger constant for a given $n$, but is bad for predicting Cheeger constant of other $n$.} on $\lambda_0(G)$ and $\lambda_1(G)$ of the spectrum and Cheeger constant data for graphs of sizes 12 and 16 for even-sized graphs and sizes 13 and 17 for odd-sized graphs. Again for training here we have taken only 40\% of the available data. Each training results in a new model, so we train the network for each $n$ a few times then take the trained model that yields the least validation error on the same $n$. We use the trained nets to predict $h(G)$ for graphs of other sizes which we compare to its true value and obtain $\Delta h$. The average deviation $\langle \Delta h \rangle$ with respect to $n$ is shown in Fig.~\ref{fig:DNNprediction} below, where we also show prediction done by linear regression method of Sec.~\ref{sec:linearpredict} for contrast. 
\begin{figure}
\centering
\includegraphics[width=70ex]{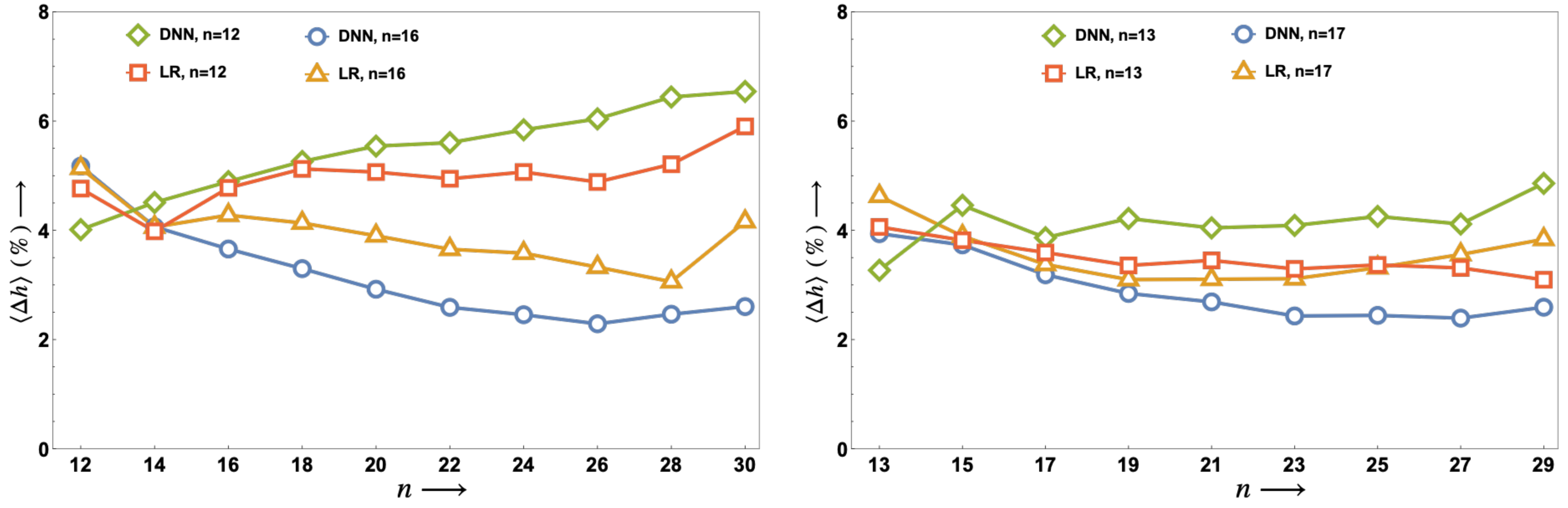}
\caption{Comparison of Deep Neural Network with Linear Regression for predicting Cheeger constant. Left panel shows mean prediction deviation for even $n$ when deep neural  network (DNN) models and linear regression (LR) models for $n=12$ or $n=16$ are used. Right panel shows mean prediction deviation for odd $n$ when deep neural  network (DNN) models and linear regression (LR) models for $n=13$ or $n=17$ are used. }
\label{fig:DNNprediction}
\end{figure}

\noindent Here are our observations
\begin{enumerate}
\item We note that $n=16$ works better than $n=12$ for predicting Cheeger constants for higher even $n$, and similarly $n=17$ works better than $n=13$ for predicting Cheeger constant for higher odd $n$.
\item Although the plots are not shown here, but we have verified that to predict for even $n$ training on even $n$ works better than training on odd $n$, and vice versa. This is consistent with observations of Sec.~\ref{sec:linearpredict}.
\item We also note that deep neural net based model provides better prediction compared to linear regression model with a consistent improvement as $n$ increases. Particularly, the models trained on $n=16$ and $n=17$ data predict Cheeger constants  for the graphs of sizes $29$ and $30$ respectively, to within 3\% accuracy on an average.
\item While we observe low average $\Delta h$ the standard deviation in $\Delta h$ is also low at about 4\% throughout the range of the $n$, thus guaranteeing reliability on predictions. This is shown in Fig.~\ref{fig:standardDeviation} below.
\end{enumerate}

\begin{figure}[H]
\centering
\includegraphics[width=55ex]{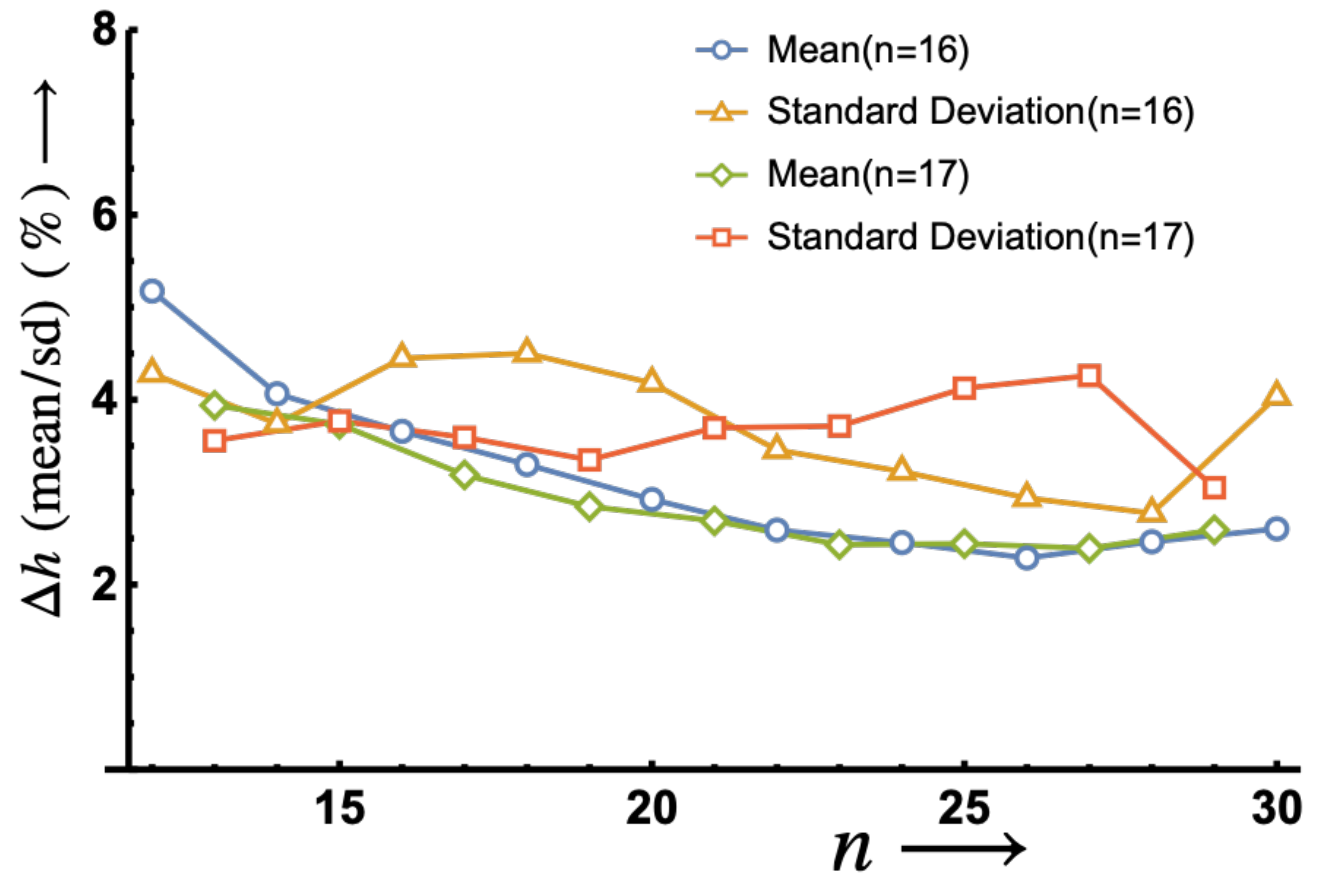}
\caption{Prediction Statistics for $n=16$ and $n=17$ DNN models. Mean deviation stays between 2\% and 4\% for all higher $n$ while standard deviation is about 4\%.}
\label{fig:standardDeviation}
\end{figure}

\section{Conclusion}
In this paper, we have studied the relevance of the spectrum of a graph $G$ in estimating $h(G)$. We find that $h(G)$ is strongly dependent on $\lambda_0(G)$ and $\lambda_1(G)$, and this correlation is largely linear with a small non-linear component, as confirmed by the machine learning analysis. We have also demonstrated that by using a deep neural network that has been moderately trained about the relationship between $h(G)$ and  $\lambda_0(G),\lambda_1(G)$, we can effectively estimate the Cheeger constant of a larger graph with high accuracy, statistically. We believe that an optimal use of this approach could be a powerful and efficient tool for studying the connectivity for large regular graphs.
\bibliographystyle{plain}
\bibliography{cheeger_estimate}

\begin{thebibliography}{10}

\bibitem{B2}
N.~Biggs.
\newblock {\em Algebraic graph theory}.
\newblock Cambridge Mathematical Library. Cambridge University Press,
  Cambridge, second edition, 1993.

\bibitem{GJS}
M.~R. Garey, D.~S. Johnson, and L.~Stockmeyer.
\newblock Some simplified {NP}-complete graph problems.
\newblock {\em Theoret. Comput. Sci.}, 1(3):237--267, 1976.

\bibitem{G1}
C.~Godsil and G.~Royle.
\newblock {\em Algebraic graph theory}, volume 207 of {\em Graduate Texts in
  Mathematics}.
\newblock Springer-Verlag, New York, 2001.

\bibitem{HLW}
S.~Hoory, N.~Linial, and A.~Wigderson.
\newblock Expander graphs and their applications.
\newblock {\em Bull. Amer. Math. Soc. (N.S.)}, 43(4):439--561, 2006.

\bibitem{KS}
Mike Krebs and Anthony Shaheen.
\newblock {\em Expander families and {C}ayley graphs}.
\newblock Oxford University Press, Oxford, 2011.
\newblock A beginner's guide.

\bibitem{LR}
Tom Leighton and Satish Rao.
\newblock Multicommodity max-flow min-cut theorems and their use in designing
  approximation algorithms.
\newblock {\em J. ACM}, 46(6):787--832, 1999.

\bibitem{LUB}
A.~Lubotzky.
\newblock {\em Discrete groups, expanding graphs and invariant measures}.
\newblock Modern Birkh\"{a}user Classics. Birkh\"{a}user Verlag, Basel, 2010.
\newblock With an appendix by Jonathan D. Rogawski, Reprint of the 1994
  edition.

\bibitem{AL}
A.~Lubotzky.
\newblock Expander graphs in pure and applied mathematics.
\newblock {\em Bull. Amer. Math. Soc. (N.S.)}, 49(1):113--162, 2012.

\bibitem{BM}
Bojan Mohar.
\newblock Isoperimetric numbers of graphs.
\newblock {\em J. Combin. Theory Ser. B}, 47(3):274--291, 1989.

\bibitem{RM1}
M.~R. Murty.
\newblock Ramanujan graphs.
\newblock {\em J. Ramanujan Math. Soc.}, 18(1):33--52, 2003.

\bibitem{SW}
A.~Steger and N.~C. Wormald.
\newblock Generating random regular graphs quickly.
\newblock volume~8, pages 377--396. 1999.
\newblock Random graphs and combinatorial structures (Oberwolfach, 1997).

\end{thebibliography}

\end{document}